\documentclass[12pt,a4paper]{article}

\usepackage{stmaryrd}
\usepackage{CJK,CJKnumb}
\usepackage{amssymb}
\usepackage{amsfonts}
\usepackage{fancyhdr}
\usepackage{bbm}
\usepackage{amsfonts}
\usepackage{amsmath}
\usepackage{amssymb}
\usepackage{indentfirst}
\usepackage{graphicx}
\usepackage{latexsym,bm, amsthm}
\usepackage{color}              
\usepackage{indentfirst}        
\usepackage{latexsym,bm}        
\usepackage{amsmath,amssymb}    
\usepackage{graphicx}
\usepackage{cases}
\usepackage{pifont}
\usepackage{txfonts}
\usepackage{fancyhdr}
\usepackage{pgf,pgfarrows,pgfnodes,pgfautomata,pgfheaps,pgfshade}
\usepackage{hyperref}
\usepackage[latin1]{inputenc}
\usepackage{colortbl}
\usepackage[english]{babel}
\usepackage{epsfig}
\usepackage{lmodern}
\usepackage[T1]{fontenc}
\usepackage{times}
\usepackage{multicol}
\usepackage{multimedia}

\newcommand{\song}{\CJKfamily{song}}            
\newcommand{\kai}{\CJKfamily{kai}}              
\newcommand{\hei}{\CJKfamily{hei}}              
\newcommand{\HZxiaowuhao}{\fontsize{9bp}{11.8bp}\selectfont}   

\allowdisplaybreaks[4]  

\CJKtilde   

\normalsize
\newcommand{\df}[2]{{\displaystyle\frac{#1}{#2}}}
\newtheoremstyle{mythm}{1.5ex plus 1ex minus .2ex}{1.5ex plus 1ex minus .2ex}{\kai}{}{\song\bfseries}{}{1em}{}
\theoremstyle{mythm}

\newcommand*{\SHOUYEJIAOZHU}[1]{\footnotetext{\hspace*{-0.7em}\HZxiaowuhao{#1}\vspace{0.4em}}} 

\begin{document}
\begin{CJK*}{GBK}{song}
\CJKindent

\title{Differential inequality of the second derivative that leads to
normality}
\author{Qiaoyu Chen, Shahar Nevo, Xuecheng Pang
}
\date{}
\maketitle \vskip 3mm

\SHOUYEJIAOZHU{2010 Mathematical  subject classification
30A10,~30D45.}
 \SHOUYEJIAOZHU{ keywords and phrases, Normal family,
Differential inequality.}

\medskip
\centerline{{\hei Abstract}}
\noindent  Let~$\mathcal{F}$~ be a
family of functions meromorphic in a domain ~D.~
 If~$\{\df{|f^{''}|}{1+|f|^{3}}:f\in \mathcal{F} \}$~is  locally
 uniformly bounded away from zero,~then~$\mathcal{F}$~ is normal.

\noindent I.~Introduction. \\
Recently,~progress was occurred concerning the
study of the connection  between differential inequalities and
normality. ~A natural point of departure for this subject is the
well-known
theorem due to F.Marty.\\
\hei{Marty's Theorem }~[8, P.75]\quad A family ~$\mathcal{F}$~of
functions meromorphic in a domain ~D~ is normal~ if and only
if~~$\{f^{\#}:~f \in \mathcal{F}\}$~is  locally
 uniformly bounded in ~D~.\\
Following Marty's Theorem,~L.~Royden proved the followiing
 generalization.\\
\hei{Theorem R}[7]\quad Let ~$\mathcal{F}$~ be a family of functions
meromorphic in a domain ~D, ~with the property that for each compact
set ~$K\subset D$~, ~there is a positive increasing function
~$h_K$~, ~such that ~$|f'(z)|\leq h_K(|f(z)|)$~for all ~$f\in
\mathcal{F} $~ ~and ~$z\in K$.~Then~$\mathcal{F}$~is normal in~D.\\
This result was significantly extended further in various
directions, ~see ~$[3],[9]~and~[11]$. ~S.Y.Li and H.Xie established
a different kind of generalization ~of Marty's Theorem that involves
higher derivatives. \\
\hei{Theorem LX} ~$[4]$\quad Let~$\mathcal{F}$~ be a family of
functions meromorphic in a domain ~D,~such that each $f\in
\mathcal{F} $~has zeros only of multiplicities $\geq k~,k\in
N$.~Then~$\mathcal{F}$~is normal in D if and only if the family
$$\left\{\df {|f^{(k)}(z)|}{1+|f(z)|^{k+1}}:f\in \mathcal{F}\right\}$$
is locally uniformly bounded in D.\\
In~$[6]$,~the second and the third authors gave a counterexample to
the validity of Theorem LX,~without the condition on the
multiplicities of zeros for ~the case ~$k=2$.\\
Concerning differential inequalities ~with the reversed sign of the
inequality,~J.~Grahl,~and the ~second author proved the ~following
result,~that may be ~considered as a counterpart to ~Marty's
Theorem.\\
\hei{Theorem GN }\quad $[1]$\quad Let ~$\mathcal{F}$~ be a family of
functions
 meromorphic in D,~and~$c>0$~.~If ~$f^{\#}(z)>c$~
 for every~$f\in \mathcal{F} $~and~$z\in D$,~then~$\mathcal{F}$~is normal in
 ~D.\\
N.Steinmetz ~$[10]$,~gave a shorter proof of Theorem GN,~using the
Schwarzian derivative and some Well-known facts on linear
differential equations.\\
Then in~$[5]$,~X.J.Liu together
 ~with the second and third
 ~authors generalized Theorem GN
 ~and proved
 the following result.\\
\hei{Theorem LNP}\quad Let~$1\leq\alpha< \infty$~ and~$c>0$.~Let
~~$\mathcal{F}$~ be the family of all meroforphic functions ~$f $
~in~D, ~such that
$${\df {|f'(z)|}{1+|f(z)|^{\alpha}}>C}$$
~for every
~$ z\in D$.\\~
Then the following hold:\\
(1) If ~$\alpha>1$,~then~$\mathcal{F}$~is normal  ~in ~D.~\\
(2) If~$\alpha=1$,~then~$\mathcal{F}$~is quasi-normal in ~D~ but
not necessarily normal. \\
Observe that (2) of the theorem is a
differential inequalities that
distinguish between quasi-normality to normality. \\
In this paper,
we continue to ~study differential inequality ~with
the reversed sign ~$(``\geq ")$~ ~and prove the following theorem.\\
\hei{Theorem 1.}\quad Let ~D~be a domain in ~$\mathbb{C}$~and let
~$c>0$.~Then the family ~~$\mathcal{F}$~ of all functions
~f~meromorphic in ~D~,~such that ~$$\df {|f^{''}(z)|}{1+|f(z)|^{3}}>
C$$~for ~every
~$z\in D$~is normal. \\
Observe that the above differential
~inequality is the reversed
inequality to that of Theorem LX in the case ~$k=2$.~\\
Let us set some notation. \\
For ~$z_0 \in C $~and ~$r>0$.~$\Delta(z_0, r) = \{z: |z - z_0| <
r\}$,~$\overline{\Delta}(z_0, r) = \{z: |z - z_0| \leq r\}$. We
write~$f_n(z)\overset\chi \Rightarrow f(z)$~on ~$D$~to indicate that
the sequence $\{f_n(z)\}$ converges to ~$f(z)$ in the spherical
metric, uniformly on compact subsets of ~$D$, and $f_n(z)\Rightarrow
f(z)$~on ~$D$ if the convergence is also in the Euclidean metric.\\
II\quad
Proof of Theorem 1.\\
Since ~$|f''| > c$~ for every~$f \in\mathcal{F},$~it ~follows
that~$\{f'':f \in\mathcal{F}\}$~is normal in ~D~.~Let
~$\{f_n\}^{\infty}_{n=1}$~be a sequence of functions from
$\mathcal{F}$.~Without loss of ~generality,~we can assume that
~$f^{''}_n(z)\overset\chi \Rightarrow H$  in D~.~Let us separate
into two cases.\\
Case 1.\quad$f_n,n\geq 1$~are holomorphic ~functions~in ~D~.\\
Case 1.1\quad H~is holomorphic function~in~D~.\\
Since normality is a local property. ~It is enough to prove that
~$\{f_n\}$~is normal at each point of ~D~. Let ~$z_0\in D$~without
loss of generality, ~we can assume that~$z_0=0$~.By the assumption
on ~H~, ~there exist some~$r >0,~M >C$,~such that~$|f^{''}_n(z)|\leq
M$~ ~for every~$z\in \Delta(0,r)$~if ~$n$~ is large ~enough.~We then
get for large enough ~$n$~ and~$z\in \Delta(0,r)$~that
~$1+|f_n(z)|^3 \leq \frac {2 M}{C}$ ~and we deduce ~that
~$\{f_n\}^{\infty}_{n=1}$~ is normal at ~$z=0,$~ ~as required.\\
Case 1.2 \quad$H \equiv\infty$~in~D.\\
Again,~let $z_0 \in D$ and assume that $z_0 =0.$~ Let $r >0$ be such
that $\overline{\Delta}(0, r)\subset D$. Without loss of
generality,~we can assume that ~$|f^{''}_n(z)|> 1$~
~for every~$z\in\Delta(0,r),~n\in\mathbb{N}$.~Then~$\log |f_n^{''}|$~is a positive harmonic function in $\Delta(0,r)$.\\
From Harnack's inequality we then get that\\
(1)~$$|f_n^{''}(z)|\leq|f_n^{''}(0)|^{\frac{1+|z|}{1-|z|}} $$~for
every $z\in\Delta(0,r),~n\in\mathbb{N}.$\\
Let us fix some~$0<\rho<\displaystyle{\frac{r}{2}}$.~Then
\\(2) $$\frac{r+\rho}{r-\rho}<3.$$\\
For every $n\geq 1,$ let ~$z_n\in\{z:|z|=\rho\}$~be such that
 $$
|f_n(z_n)|=\max\limits_{|z|\leq\rho}|f_n(z)|=M(\rho,f_n)$$
 By Cauchy's Inequality ,~we get that\\
$$|f^{''}_n(0)|\leq
\df{2}{\rho^2}M(\rho,f_n)=\df{2}{\rho^2}|f_n(z_n)|.$$
 Hence,~by (1),~we get
$$C\leq \df{|f^{''}_n(z_n)|}{1+|f_n(z_n)|^3}\leq
\df{|f^{''}_n(z_n)|}{|f_n(z_n)|^3}\leq
\df{|f^{''}_n(0)|^{\df{r+\rho}{r-\rho}}}{|f_n(z_n)|^3} \leq
 \left(\df{2}{\rho^2}\right)^{\df{r+\rho}{r-\rho}}\left|f_n(z_n)\right|^{\df{r+\rho}{r-\rho}-\displaystyle{3}},$$
Thus, by (2)
$$M(\rho,f_n)=|f(z_n)|\leq \left(\df{1}{C}\left(\df{2}{\rho^2}\right)^{\df{r+\rho}
{r-\rho}}\right)^{\df{1}{3-\df{r+\rho}{r-\rho}}},$$ which means that
$\{f_n\}$ is locally uniformly bounded in $\Delta(0,\rho)$ and thus
$\{f_n\}$ is normal at $z=0$.\\
Case 2 \quad $f_n$ are meromorphic functions with pole in $D$.\\
By Case 1 we have to prove normality only at point $z_0$,~where
$H(z_0)=\infty$.~Such points exist if $H$ is a meromorphic function
with poles in $D$ or if $H\equiv \infty$. So let $z_0$ be such that
$H(z_0)\equiv \infty$.~Without loss of generality , we can assume
that $z_0=0$. After moving to a subsequence,  that without loss of
generality will also be denoted by $\{f_n\}_1^{\infty}$,~we can
assume that there is a sequence $\zeta_n\rightarrow 0$~such that
~$f_n(\zeta_n)=\infty$.~For if it was not the case,then for some
$\delta
>0$ ~and large enough ~$n$~,$f_n$~would be holomorphic in
$\Delta(0,\delta)$,and then we would get the asserted normality by
case (1).\\
Also we can assume the existence of \\
(3)\quad\quad a sequence~$\eta_n\rightarrow 0$~such that $f_n(\eta_n)=0$.\\
Indeed,~since ~$H(z_0)=\infty$~there exists some ~$\delta> 0$~such
that for large enough ~$n$~$\min\limits_{z\in \Delta
(0,\delta)}|f_n^{''}|>1$.\\
Combining it with ~$f_n\neq 0$~ in some neighbourhood of
~$z=0$~gives the normality at~$z=0$~by Gu's Criterion [2]. \\
We can also assume that ~$\{f_n^{'}\}$~ is not normal
at~$z=0$.~Indeed, if~$\{f_n^{'}\}$~would be normal at ~$z=0$,~then
by Marty's theorem there exist ~$r_1>0$~and ~$M>0$~such that for
large enough ~$n$,~ $\df{|f_n^{''}(z)|}{1+|f_n^{'}(z)|^2}<M  $ for
~$z\in \Delta(0,r_1)$.~ Since ~$H(0)=\infty $,there exists some
~$r_2\leq r_1$~such that for large enough ~$n,$~$|f_n^{''}(z)|\geq
2M$~for
~$z\in\Delta(0,r_2)$.\\
We thus have for large enough ~$n$~and ~$z\in \Delta(0,r_2)$,~
$1+|f_n^{'}(z)|^2>\df{|f_n^{''}(z)|}{M}\geq 2$~and thus
~$|f^{'}_n(z)|\geq 1$.We then get
$$\df{|f^{'}_n(z)|^2}{|f^{''}_n(z)|}=\df{|f^{'}_n(z)|^2}{1+|f^{'}_n(z)|^2}\cdot
\df{1+|f^{'}_n(z)|^2}{|f^{''}_n(z)|}\geq \df{1^2}{1+1^2}\cdot
\df{1}{M}=\df{1}{2M}.$$ Hence We have for large enough ~$n$~and ~
$z\in
\Delta(0,r_2)$\\
$(4)\quad\quad\displaystyle{\df{|f^{'}_n(z)|^2}{1+|f_n(z)|^3}=\df{|f^{'}_n(z)|^2}{|f^{''}_n(z)|}\cdot
\df{|f^{''}_n(z)|}{1+|f_n(z)|^3}>\df{1}{2M}\cdot C}.$ \\
Now,for every ~$x\geq 0,\df{\sqrt{1+x^2}}{1+x}\geq
\df{1}{\sqrt{2}}$,~and by taking square root of (4),~we get
$$\df{|f^{'}_n(z)|}{1+|f_n(z)|^\frac{3}{2}}=\df{|f^{'}_n(z)|}{\sqrt{1+|f_n(z)|^3}}\cdot
\df{\sqrt{1+|f_n(z)|^3}}{1+|f_n(z)|^\frac{3}{2}}>\sqrt{\df{C}{2M}}\cdot\df{1}{\sqrt{2}}
.$$ \\
By (1) of Theorem LNP, ~with ~$\alpha=\frac{3}{2}>1$,~we
deduce that~$\{f_n\}$ ~ is normal in $\Delta(0,r_2)$ and~we are
done.
\\
Thus
we can assume that~$\{f'_n\}$ ~is not normal at ~$z=0$. \\
Similarly
to (3) we can assume~that there is a sequence
~$s_n\to 0$~such that~$f^{'}_n(s_n)=0$.\\
We claim that we
can  assume that~$\{\df{f^{'}_n}{f^{''}_n}\}_{n=1}^{\infty}$~is not
normal at
~$z=0$. \\
Otherwise,~after moving to a~subsequence that will also be denoted
by~$\{\df{f^{'}_n}{f^{''}_n}\}_{n=1}^{\infty}$ we
have~$\df{f^{'}_n}{f^{''}_n}\Rightarrow H_1$ in $\Delta(0,r)$ , for
some~$r
> 0$.~Since ~$f^{''}_n \neq 0$~and
$\df{f^{'}_n}{f^{''}_n}(\zeta_n)=0$~ then ~$H_1$~must be holomorphic
function in~$\Delta(0,r).$
Differentiation then gives\\
(5)\quad\quad $1-\df{f^{'}_n f^{''}_n }{(f^{''}_n)^2}\Rightarrow
H^{'}_1$~in $\Delta(0,r).$\\
At ~$z=s_n$~ the left hand of~(5)~is equal to 1.~on the other hand
~in some small neighbourhood of~$z=\zeta_n$,~We have
$f_n(z)=\frac{A_n}{z-z_n}+\hat{f}_n(z),~$ where $A_n\neq 0$ is a
constant,~and $\hat{f}_n(z)$ is analytic. ~Here we used that
according to the ~assumption of Theorem ~$1$~,~all ~poles of
~$f_n$~must be simple. \\
Hence we have
~$f^{'}_n(z)=\df{-A_n}{(z-\zeta_n)^2}+\hat{f}_n^{'}(z),
f^{''}_n(z)=\df{2A_n}{(z-\zeta_n)^3}+\hat{f}_n^{''}(z),
f_n^{(3)}(z)=\df{-6A_n}{(z-z_n)^4}+\hat{f}_n^{(3)}(z)$.~ Then the
left hand of (5) get at~$z=\zeta_n.$~The value
~$1-\df{6}{4}=-\df{1}{2}\neq 1$,
a contradiction .\\
Claim ~ there exist ~$r>0$~and ~$k>0$~such that for large enough
~$n$,
$|\df{f_n}{f_n^{''}}(z)|,~~\left|\df{f_n^2}{f_n^{''}}(z)\right|\leq K$~for ~$z\in \Delta(0,r)$.\\
Proof of Claim \quad Since ~$H(0)=\infty $,there exist ~$r>0$~and
~$M>0$~such that $\overline{\Delta}(0,r)\subset D$~and such that for
large enough $n$, $|f_n^{''}(z)|>M$ for ~$z\in \Delta (0,r)$.\\
Now ,if $|f_n(z)|\leq |f_n^{''}(z)|^{\frac{1}{3}}$ then \\
(6) \quad\quad $|\df{f_n}{f_n^{''}}(z)|\leq
\df{|f_n^{''}(z)|^{\frac{1}{3}}}{|f_n^{''}(z)|}\leq
\df{1}{M^{\frac{2}{3}}}$\\
and\\
(7)\quad\quad$|\df{f_n^2}{f_n^{''}}(z)|\leq
\df{|f_n^{''}(z)|^{\frac{2}{3}}}{|f_n^{''}(z)|}\leq
\frac{1}{M^{\frac{1}{3}}}$.\\
If on the other hand $|f_n(z)|\geq
|f_n^{''}(z)|^{\frac{1}{3}}$,~then since $\df{x}{1+x^3}\leq
\df{2^{\frac{2}{3}}}{3}$~for ~$x\geq 0$,~we
get \\
(8)\quad\quad
$|\df{f_n}{f_n^{''}}(z)|=\df{1+|f_n(z)|^3}{|f_n^{''}(z)|}\cdot
\df{|f_n(z)|}{1+|f_n(z)|^3}\leq \df{1}{C}\cdot
\df{2^{\frac{2}{3}}}{3}$.\\
Also We have $\df{x^2}{1+x^3}\leq \frac{2^{\frac{2}{3}}}{3}$~for
~$x\geq
0$ and thus\\
(9)\quad\quad$|\df{f_n^2}{f_n^{''}}(z)|=\df{1+|f_n(z)|^3}{|f_n^{''}(z)|}\cdot
\df{|f_n^2(z)|}{1+|f_n(z)|^3}\leq \df{1}{C}\cdot
\df{2^{\frac{2}{3}}}{3}$.\\
The claim then follows by taking
$k=\max{\{\frac{1}{M^{\frac{2}{3}}},\frac{1}{M^{\frac{1}{3}}},\frac{1}{C}\cdot\frac{2^{\frac{2}{3}}}{3}\}}$
and consider (6),(7),(8)and (9).\\
From the claim we deduce that $\{\df{f_n}{f_n^{''}}\}_{=1}^\infty$
and $\{\df{f_n^2}{f_n^{''}}\}_{=1}^\infty$ are normal in
$\Delta(0,r)$,~so after moving to a subsequence,~that also will be
denote by
$\{f_n\}_{n=1}^\infty $,~we get that \\
(10)\quad\quad$\df{f_n}{f_n^{''}}\rightarrow H_1$ in $\Delta(0,r)$\\
and\\
(11)\quad\quad$\df{f_n^2}{f_n^{''}}\rightarrow H_2$ in $\Delta(0,r)$\\
From the claim it follows that $H_1$ and $H_2$ are holomorphic in
$\Delta(0,r)$.\\
Differentiating (10) and (11) gives respectively \\
(12)\quad\quad$\df{f_n^{'}}{f_n^{''}}-\df{f_n^{(3)}}{f_n^{''}}\cdot
f_n
\Rightarrow H_1^{'}$ in $\Delta(0,r)$\\
and \\
(13)\quad\quad$2f_n\cdot\df{f_n^{'}}{f_n^{''}}-f_n^2\cdot\df{f_n^{(3)}}{{f^{''}_n}^2}\Rightarrow
H_2^{'}$ in $\Delta(0,r)$.\\
Since  $\{f_n''\}_{n=1}^\infty$ is normal, there exists some $k_1>0$
such that $ \frac{|f_n^{(3)}(z)|}{1+|f_n''(z)|}\leq k_1 $ for every
$n\geq1$ and for every $z\in\Delta(0,r)$. Since in addition for
large enough $n$, $|f_n''(z)|>M$, then
\begin{eqnarray*}
\frac{|f_n^{(3)}(z)|}{|f_n^{''}(z)|^2}&=&\frac{|f_n^{(3)}(z)|}{1+|f_n^{''}(z)|^2}\frac{1+|f_n^{''}(z)|^2}{|f_n^{''}(z)|^2}\\
&\leq&k_1(1+\frac{1}{M^2}):=k_2.
\end{eqnarray*}
\\
Thus~\\
(14)\quad\quad$
\frac{|f_n^{''}(z)|^2}{|f_n^{(3)}(z)|}\geq\frac{1}{k^2} $~for large
enough n. \\
Now since we assume that $\{\frac{f_n'}{f_n''}\}$ is not normal at
$z=0$, then after moving to a subsequence, that also will be denoted
by $\{f_n\}_{n=1}^\infty$, we get that there exists a sequence of
points $t_n\rightarrow0$, such that
$$
\frac{f_n'}{f_n''}(t_n):=M_n\rightarrow\infty, \quad M_n\in
\mathbb{C}.
$$
Substituting $z=t_ n$ in $(12)$ gives\\
 (15)
$$
M_n-\frac{f_n^{(3)}\cdot
f_n}{f_n''^2}(t_n):=\varepsilon_n\rightarrow H_1'(0).
$$
Hence
$$
f_n(t_n)=(M_n-\varepsilon_n)\frac{f_n^{''2}}{f_n^{(3)}}(t_n)
$$
From $(15)$ we get,~ by substituting $z=t_n$ in $(13)$
$$
2(M_n-\varepsilon_n)\frac{f_n^{''2}}{f_n^{(3)}}(t_n)M_n-(M_n-\varepsilon_n)^2\left(\frac{f_n^{''2}}{f_n^{(3)}}(t_n)\right)^2\frac{f_n^{(3)}}{f_n^{''2}}(t_n):=\delta_n\rightarrow
H_2'(0).
$$
\\
From this we get after simplifying
$$
(M_n^2-\varepsilon_n^2)\frac{f_n^{''2}}{f_n^{(3)}}(t_n)=\delta_n.
$$
But by (14) the left hand above tends to $\infty$~ as
$n\rightarrow\infty$, while the right hand is bounded, a
contradiction.\\
This completes the proof of Theorem 1.

QIAOYU CHEN, DEPARTMENT OF MATHEMATICS, EAST CHINA NORMAL
UNIVERSITY,SHANG HAI 200241,P.R.CHINA\\
E-mail address: goodluckqiaoyu@126.com

SHAHAR NEVO, DEPARTMENT OF MATHEMATICS, BAR-ILAN UNIVERSITY, 52900
RAMAT-GAN, ISRAEL \\
E-mail address: nevosh@macs.biu.ac.il

XUECHENG PANG, DEPARTMENT OF MATHEMATICS, EAST CHINA NORMAL
UNIVERSITY,SHANG HAI 200241,P.R.CHINA\\
E-mail address: xcpang@math.ecnu.cn

\clearpage

\end{CJK*}
\end{document}